\documentclass[11pt,reqno]{amsart}
\usepackage{mathrsfs}
\usepackage{amsmath}
\usepackage{amssymb, amsmath}
\usepackage{graphicx}

\newcommand{\newcom}{\newcommand}

\newcom{\al}{\alpha}
\newcom{\be}{\beta}
\newcom{\eps}{\epsilon}
\newcom{\ga}{\gamma}
\newcom{\Ga}{\Gamma}
\newcom{\ka}{\kappa}
\newcom{\Lam}{\Lambda}
\newcom{\lam}{\lambda}
\newcom{\Om}{\Omega}
\newcom{\om}{\omega}
\newcom{\Si}{\Sigma}
\newcom{\si}{\sigma}
\newcom{\tht}{\theta}
\newcom{\dtri}{\nabla}
\newcom{\tri}{\triangle}
\newcom{\oo}{\infty}
\newcom{\vphi}{\varphi}
\newcom{\cB}{{\mathcal B}}
\newcom{\cC}{{\mathcal C}}
\newcom{\cD}{{\mathcal D}}
\newcom{\cF}{{\mathcal F}}
\newcom{\cL}{{\mathcal L}}
\newcom{\cM}{{\mathcal M}}
\newcom{\cP}{{\mathcal P}}
\newcom{\cS}{{\mathcal S}}
\newcom{\cQ}{{\mathcal Q}}
\newcom{\caly}{{\mathcal Y}}
\newcom{\calz}{{\mathcal Z}}
\newcom{\bfz}{{\bf Z}}
\newcom{\R}{\Bbb R}
\newcom{\N}{\Bbb N}
\newcom{\Z}{\Bbb Z}
\newcom{\C}{\Bbb C}
\newcom{\E}{\Bbb E}

\newcom{\f}{\frac}
\newcom{\di}{\displaystyle\int}
\newcom{\ds}{\displaystyle\sum}
\newcom{\dl}{\displaystyle\lim}
\newcom{\ov}{\overline}
\newcom{\sset}{\subset}
\newcom{\wt}{\widetilde}
\newcom{\pa}{\partial}
\newcom{\p}{\partial}
\newcom\na{\nabla}
\newcom{\co}{\cdot}
\newcom{\suml}{\sum\limits}
\newcom{\supl}{\sup\limits}
\newcom{\intl}{\int\limits}
\newcom{\infl}{\inf\limits}
\newcom{\disp}{\displaystyle}
\newcom{\non}{\nonumber}
\newcom{\no}{\noindent}
\newcom{\QED}{$\square$}

\def\ef{\hphantom{MM}\hfill\llap{$\square$}\goodbreak}

\def\eqdefa{\buildrel\hbox{\footnotesize def}\over =}

\newtheorem{athm}{\bf \t}[section]
\newenvironment{thm} [1] {\def\t{#1}\begin{athm} \bf \rm} {\end{athm}}
\newcom{\bthm}{\begin{thm}}\newcom{\ethm}{\end{thm}}

\newcom{\beq}{\begin{equation}}
\newcom{\eeq}{\end{equation}}
\newcom{\ben}{\begin{eqnarray}}
\newcom{\een}{\end{eqnarray}}
\newcom{\beno}{\begin{eqnarray*}}
\newcom{\eeno}{\end{eqnarray*}}
\newcom{\bali}{\begin{aligned}}
\newcom{\eali}{\end{aligned}}

\numberwithin{equation}{section}

\begin{document}

\title[uniqueness for the Navier-Stokes equations]
{on the uniqueness of weak solutions for the 3D Navier-Stokes
equations}

\author{Qionglei Chen}
\address{Institute of Applied Physics and Computational Mathematics,P.O. Box 8009, Beijing 100088, P. R. China}
\email{chen\_qionglei@iapcm.ac.cn}

\author{Changxing Miao}
\address{Institute of Applied Physics and Computational Mathematics,P.O. Box 8009, Beijing 100088, P. R. China}
\email{miao\_changxing@iapcm.ac.cn}

\author{Zhifei Zhang}
\address{School of Mathematical Sciences, Peking University, 100871, P. R. China}
\email{zfzhang@math.pku.edu.cn}


\date{22,August}



\begin{abstract}
In this paper, we improve some known uniqueness results of weak
solutions for the 3D Navier-Stokes equations. The proof uses the
Fourier localization technique and the losing derivative estimates.
\end{abstract}

\maketitle

\section{Introduction}
 We consider the three dimensional Navier-Stokes equations in $\R^3$
\begin{equation}\label{NSequ}
\left\{
\begin{array}{ll}
u_t-\Delta u+u\cdot\nabla u+\nabla p=0,\\
\textrm{div}\, u = 0,\\
u(0)=u_0(x),
\end{array}
\right.
\end{equation}
where $u=(u^1(t,x),u^2(t,x),u^3(t,x))$ and $p=p(t,x)$ denote the
unknown velocity vector and  the unknown scalar pressure of the
fluid  respectively, while
$u_0(x)$ is a given initial velocity vector satisfying $\textrm{div}
u_0 = 0$.

In a seminal paper \cite{Ler}, J. Leray  proved the global existence
of weak solution with finite energy, that is,
$$u(t,x)\in \cL_T\eqdefa L^{\infty}(0,T;L^2)\cap L^2(0,T;H^1)\quad  \textrm{for any }T>0.$$
It is well known that weak solution is unique and regular in two
spatial dimensions. In three dimensions, however, the question of
regularity and uniqueness of weak solution is an outstanding open
problem in mathematical fluid mechanics. In this paper, we are
interested in the classical problem of finding sufficient conditions
for weak solutions of (\ref{NSequ}) such that they become regular
and unique. Let us firstly recall the definition of weak solution.

\bthm{Definition} Let $u_0 \in L^2(\R^3)$ with
$\textrm{div}u_0=0$. A measurable  function $u$ is called a weak
solution of (\ref{NSequ}) on $(0,T)\times \R^3$ if it satisfies
the following conditions:
\begin{flushleft}
\textrm{(1)}\, $u\in \cL_T\cap C_w([0,T];L^2)$, where
$C_w([0,T];L^2)$ consists of  all weak continuous functions with
respect to time in $L^2(\R^3)$;

\textrm{(2)}\, $\textrm{div}\,u=0$ in the sense of distribution;

\textrm{(3)}\, For any function $\psi \in C^\infty_0([s,t]\times
\R^3)$ with $\textrm{div}\, \psi=0$, there holds \beno
&&\int_s^t\int_{\R^3}\{u\cdot \psi_t-\na u\cdot \na \psi+\na \psi:(u\otimes u)\}(t',x)dxdt'\\
&&\qquad=\int_{\R^3}u(t,x)\cdot\psi(t,x)dx-\int_{\R^3}u(s,x)\cdot\psi(s,x)dx;
\eeno
\end{flushleft}
In addition, if $u$ satisfies the energy inequality \beno
\|u(t)\|_2^2+2\int_0^t\|\na u(t')\|_2^2dt'\le \|u_0\|_2^2,\quad
\eeno it is also called a Leray-Hopf weak solution. \ethm

The Leray-Hopf weak solutions are unique and regular in the class
\beno &&{\cP}=L^q(0,T;L^r) \quad \textrm{with} \quad \frac 2 q+\frac
3 r=
1,\, 3\le r\le \infty,\quad [11,14,15,25,27]  \\
&& \textrm{or}\quad {\cP}=L^q(0,T;W^{1,r}) \quad \textrm{with} \quad
\frac 2 q+\frac 3 r=
2,\, \f32<r\le \infty,\quad[1] \\
&& \textrm{or}\quad {\cP}= L^q(0,T;W^{s,r}) \quad \textrm{with}
\quad \frac 2 q+\frac 3 r= 1+s,\, \f3{1+s}<r\le \infty, s\ge
0,\,[26]. \eeno Recently, there are many researches devoted to
refine the above results. First of all, we have the following
refined regularity criterion in the framework of Besov spaces: the
weak solutions are regular in the class \beno {\cP}=
C([0,T];B^{-1}_{\infty,\infty}) \quad \textrm{or}\quad  {\cP}=
L^q(0,T; B^r_{p,\infty}), \eeno with $\f 2 q+\f 3 p =1+r, \f 3
{1+r}<p\le \infty,$ and $-1<r\le 1$, see \cite{CZ,Che,KOT,KSI}.
Concerning the refined uniqueness criterion of weak solutions,
Kozono and Taniuchi\cite{KT} proved the uniqueness of the
Leray-Hopf weak solutions in the class \beno {\cP}= L^2(0,T;BMO);
\eeno Gallagher and Planchon\cite{GP} proved the uniqueness in the
class \beno {\cP}= L^q(0,T; \dot B^{-1+\f3p+\f2q}_{p,q})\quad
\textrm{with}\quad \f3p+\f2q>1; \eeno
Lemari\'{e}-Rieusset\cite{Lem1} proved the uniqueness in the class
\beno {\cP}= C([0,T];X_1^{(0)})\quad \textrm{or}\quad {\cP}=
L^{\f2{1-r}}(0,T; X_r)\quad \textrm{with}\quad r\in [0,1); \eeno
Finally, Germain\cite{Ger} proved the uniqueness  in the class
\beno {\cP}= C([0,T];X_1^{(0)})\quad \textrm{or}\quad  {\cP}=
L^{\f2{1-r}}(0,T; X_r)\quad \textrm{with}\quad r\in [-1,1). \eeno
Here $B^s_{p,q}$ denotes the Besov space and \beno X_s:= \left\{
\begin{array}{ll}
\textbf{M}(\dot H^s,L^2),\quad \textrm{if}\,\, s\in (0,1]\\
\Lambda^sBMO,\quad \textrm{if}\,\, s\in (-1,0]\\
\textbf{Lip},\quad \textrm{if}\,\, s=-1,
\end{array}
\right. \eeno where $\textbf{M}(\dot H^s,L^2)$ is the space of
distributions such that their pointwise product with a function in
$\dot H^s$ belongs to $L^2$, $\Lambda^s=(1-\Delta)^{\f s2}$.
$X_s^{(0)}$ denotes the closure of the Schwartz class in $X_s$. We
want to point out that \ben\label{embedding}
 X_s\hookrightarrow \Lambda^s BMO,\quad
\textrm{if}\,\, s\in (0,1]. \een We refer to \cite{Ger} for more
properties about $X_s$. The key step of their proofs is to find a path space  $\cP$ so that
the trilinear form \beno F(u,v,w):=\int_0^T\int_{\R^3}u\cdot \na
v\cdot w dxdt \eeno is continuous from $(\cL_T)^2\times {\cP}$ to
$\R$. Germain also pointed out that the path space $\cP$ he found is
optimal in some sense( see P.400 in \cite{Ger} for precise meaning).

The purpose of this paper is to improve the above uniqueness results.

\bthm{Theorem}\label{thm1} Let $u_0,v_0
\in L^2(\R^3)$ with $\textrm{div}u_0=\textrm{div} v_0=0$. Let $u$
and $v$ be two Leray-Hopf weak solutions of
 (\ref{NSequ}) on $(0,T)$ with the initial data $u_0$ and $v_0$ respectively. Assume that
\beno u\in L^{q}(0,T; B^{r}_{p,\infty}), \eeno with $\f
2{q}+\f3{p}=1+r, \f 3 {1+r}<p\le \infty, r\in (0,1]$, and $(p,r)\neq
(\infty,1)$. Then there holds \beno
&&\|u(t)-v(t)\|_2^2+\int_0^t\|\na (u-v)(t')\|_2^2dt'\\
&&\qquad\le
\|u_0-v_0\|_2^2\exp\bigg\{C\int_0^t(e+\|u(t')\|_{B^{r}_{p,\infty}})^qdt'\bigg\}.
\eeno In particular, if $u_0=v_0$, then $u=v$ a.e. on $(0,T)\times
\R^3$. \ethm

\bthm{Remark} Due
to the embedding relation \beno B^s_{p,q}\subsetneq
B^s_{p,\infty}\quad q<+\infty \quad \textrm{and}
 \quad \Lambda^{-r}BMO\subsetneq B^{r}_{\infty,\infty},
\eeno Theorem \ref{thm1} is an improvement of the corresponding
results given by  Gallagher and Planchon \cite{GP} and  Germain
\cite{Ger}.  The proof only uses an important observation that if
$u\in L^{q}(0,T; B^{r}_{p,\infty})$ with $(p,q,r)$ as in Theorem
\ref{thm1}, then $u$ can be decomposed as
$$
u=u^l+u^h\quad \textrm{with }u^l\in L^1(0,T; \textbf{Lip}) \textrm{
and } u^h\in L^{\widetilde{q}}(0,T; L^{\widetilde{p}})
$$
for some $\widetilde{}p,\widetilde{q}$ satisfing $\f 2
{\widetilde{q}}+\f3{\widetilde{p}}=1,\widetilde{p}>3$, see Lemma \ref{decom}. \ethm

In the case of $r\le0$ and $(p,r)=(\infty,1)$, using Bony's
decomposition and the losing derivative estimates, we prove

\bthm{Theorem}\label{thm2} Let $u_0 \in L^2(\R^3)$ with $\textrm{div}u_0=0$. Let $u$ and $v$ be two weak solutions of
 (\ref{NSequ}) on $(0,T)$ with the same initial data $u_0$. Assume that $u$ and $v$
 satisfy one of the following two conditions:

(a) $u\in L^{q_1}(0,T; B^{r_1}_{p_1,\infty})$ and $v\in L^{q_2}(0,T;
B^{r_2}_{p_2,\infty})$, where
 $$\f 2{q_1}+\f3{p_1}=1+r_1,\quad\f 2{q_2}+\f3{p_2}=1+r_2,$$
with $r_1, r_2\in (-1,0], r_1+r_2>-1, \f 3
{1+r_1}<p_1\le \infty, \f 3 {1+r_2}<p_2\le \infty$.

(b) $u,v\in L^{1}(0,T; B^{1}_{\infty,\infty})$.

\no Then $u=v$ a.e. on $(0,T)\times \R^3$. \ethm

\bthm{Remark} Due to the embedding relation
$$
X_s\subsetneq \Lambda^s BMO\subsetneq B^{-s}_{\infty,\infty}, \quad
s\in (0,1],
$$
the condition imposed on weak solution in Theorem \ref{thm2} is
weaker than that of Germain \cite{Ger} and
Lemari\'{e}-Rieusset\cite{Lem1}. However, the price to pay is to
impose the conditions on both weak solutions. \ethm

\bthm{Remark} The main novelty of Theorem \ref{thm2} is that weak
solutions are uniqueness in the class $L^{1}(0,T;
B^{1}_{\infty,\infty})$. In particular, from the inequality
\ben\label{BKM-inequ} \|u\|_{B^1_{\infty,\infty}}\le
C(\|u\|_2+\|\textrm{curl} u\|_{B^0_{\infty,\infty}}),\quad
(\textrm{see Section 2 for its proof}) \een we can obtain the
\textbf{Beale-Kato-Majda} type uniqueness criterion: if weak
solutions $u$ and $v$ with the same initial data satisfy
$$
\textrm{curl} u, \textrm{curl} v\in L^{1}(0,T;
B^{0}_{\infty,\infty}),
$$
then $u=v$ on $(0,T)\times \R^3$. Secondly, Theorem \ref{thm2} allows us to
impose different conditions on both weak solutions. Thirdly, we
don't impose the energy inequality on weak solutions.

\ethm

\bthm{Remark} Chemin and Lemari\'{e}-Rieusset\cite{Ch1, Lem2}
proved the uniqueness of weak solutions in the class
$C([0,T];B^{-1}_{\infty,\infty})$. While, Theorem \ref{thm2} gives
the uniqueness in the class $L^1(0,T;B^{1}_{\infty,\infty})$. It
is natural to expect that the uniqueness also holds in the class
$L^{\f2{1+r}}(0,T;B^{r}_{\infty,\infty})$ for $r\in (-1,1)$ from
the viewpoint of interpolation. This problem remains unknown for
the case of $r\in (-1,-\f12]$. \ethm

\noindent{\bf Notation.} Throughout the paper, $C$ stands for   a
generic  constant. We  will use the notation $A\lesssim B$ to denote
the relation  $A\le CB$, and
 $\|\cdot\|_{p}$ denotes the norm of the Lebesgue space $L^p$.

\section{Preliminaries}

Let us firstly recall some basic facts on the Littlewood-Paley
decomposition, one may check \cite{Ch} for more details. Choose
two nonnegative radial functions $\chi$, $\varphi \in {\cS}(\R^3)$
supported respectively in ${\cB}=\{\xi\in\R^3,\,
|\xi|\le\frac{4}{3}\}$ and ${\cC}=\{\xi\in\R^3,\,
\frac{3}{4}\le|\xi|\le\frac{8}{3}\}$ such that for any $\xi\in
\R^3$, \ben\label{unitydecompositon}
\chi(\xi)+\sum_{j\ge0}\varphi(2^{-j}\xi)=1. \een Let
$h={\cF}^{-1}\varphi$ and $\tilde{h}={\cF}^{-1}\chi$, the
frequency localization operator $\Delta_j$ and $S_j$ are defined
by
\begin{align}
&\Delta_jf=\varphi(2^{-j}D)f=2^{3j}\int_{\R^3}h(2^jy)f(x-y)dy,\quad\mbox{for}\quad
j\geq 0, \nonumber\\
&S_jf=\chi(2^{-j}D)f=\sum_{-1\leq k\le
j-1}\Delta_kf=2^{3j}\int_{\R^3}\tilde{h}(2^jy)f(x-y)dy,\quad\mbox{and}\nonumber\\
&\Delta_{-1}f=S_{0}f, \qquad\Delta_{j}f=0\quad\mbox{for}\quad j\le
-2. \nonumber\end{align}
With our choice of $\varphi$, one can easily verify that  \beq\label{orthproperty}
\begin{aligned}
&\Delta_j\Delta_kf=0\quad \textrm{if}\quad|j-k|\ge 2\quad
\textrm{and}
\quad \\
&\Delta_j(S_{k-1}f\Delta_k f)=0\quad \textrm{if}\quad|j-k|\ge 5.
\end{aligned}
\eeq For any $f\in \cS'(\R^3)$, we have by (\ref{unitydecompositon}) that
\ben\label{Littlewoodpaleydeom} f=S_0(f)+\sum_{j\ge 0}\Delta_jf,
\een which is called the Littlewood-Paley decomposition. In the
sequel, we will constantly use the Bony's decomposition from
\cite{Bon} that \beq\label{Bonydecom} uv=T_uv+T_vu+R(u,v), \eeq with
$$T_uv=\sum_{j}S_{j-1}u\Delta_jv, \quad R(u,v)=\sum_{|j'-j|\le 1}\Delta_ju \Delta_{j'}v,$$
and we also denote
$$T'_{u}v=T_{u}v+R(u,v).$$

With the introduction of $\Delta_j$, let us recall the definition of
the inhomogenous Besov space from \cite{Tri}: \bthm{Definition} Let
$s\in \R, 1\le p,q\le\infty$, the inhomogenous Besov space
$B^s_{p,q}$ is defined by
$$B^s_{p,q}=\{f\in {\cS}'(\R^3); \quad \|f\|_{B^s_{p,q}}<\infty\},$$ where
$$\|f\|_{B^s_{p,q}}:=\left\{\begin{array}{l}
\displaystyle\bigg(\sum_{j= -1}^{\infty}2^{jsq}\|\Delta_j
f\|_{p}^q\bigg)^{\frac 1
q},\quad \hbox{for}\quad q<\infty,\\
\displaystyle\sup_{j\geq -1}2^{js}\|\Delta_jf\|_{p}, \quad \hbox{
for} \quad q=\infty.
\end{array}\right.
$$
\ethm Let us point out that $B^{s}_{\infty,\infty}$ is the usual
H\"{o}lder space $C^{s}$ for $s\in \R\setminus\Z$ and the following
inclusion relations hold \beno \textbf{Lip}\subsetneq
B^{1}_{\infty,\infty}, \quad \Lambda^{-s}BMO\subsetneq
B^{s}_{\infty,\infty} \quad \textrm{for } s\in \R. \eeno We refer to
\cite{Ger,Tri} for more properties.

The following Bernstein's inequalities will be frequently used
throughout the paper.

\bthm{Lemma}\cite{Ch}\label{Berstein} Let $1\le p\le q\le\infty$.
Assume that $f\in L^p$, then there hold
\beno &&{\rm supp}\hat{f}\subset\{|\xi|\leq
C2^j\}\Rightarrow \|\partial^\alpha f\|_{q}\le
C2^{j{|\alpha|}+3j(\frac{1}{p}-\frac{1}{q})}\|f\|_{p},
\\
&&{\rm supp}\hat{f}\subset\{\f1{C}2^j\leq |\xi|\leq C
2^j\}\Rightarrow \|f\|_{p}\le
C2^{-j|\alpha|}\sup_{|\beta|=|\al|}\|\partial^\beta f\|_{p}. \eeno
Here the constant $C$ is independent of $f$ and $j$.
\ethm

We conclude this section by a proof of the inequality
(\ref{BKM-inequ}). Using  Lemma \ref{Berstein}, we have \beno
\|\Delta_{-1}u\|_{\infty}\le C\|u\|_2, \eeno and  for $j\ge 0$,
\beno 2^j\|\Delta_{j}u\|_{\infty}\le C\|\Delta_j\na u\|_\infty.
\eeno Due to the Biot-Savart law\cite{Maj}, $\nabla u$ can be
written as \beno \nabla u(x)=Cw(x)+K\ast w(x),\quad
w=\textrm{curl}\, u \eeno where $C$ is a constant matrix, and $K$
is a matrix valued function with homogeneous of degree $-3$. So,
we get that for $j\ge 0$, \beno 2^j\|\Delta_{j}u\|_{\infty}\le
C\|\Delta_j w\|_\infty, \eeno where we used the fact that
$$
\|\Delta_j(Tf)\|_p\le C\|\Delta_jf\|_p,\quad \textrm{for} \quad j\ge
0, 1\le p\le \infty,
$$
if $T$ is a  singular integral operator of convolution type with
smooth kernel \cite{Ste}. Then the inequality
(\ref{BKM-inequ}) is concluded from the definition of Besov space.

\section{Proofs of Theorems}

 This section is devoted to the proof of Theorem \ref{thm1} and Theorem \ref{thm2}.

\subsection{The proof of Theorem \ref{thm1}}

The proof is based on the following
decomposition lemma which may be independent of interest.
\bthm{Lemma}\label{decom}Assume that $u\in L^{q}(0,T;
B^{r}_{p,\infty})$ with $\f 2{q}+\f3{p}=1+r, \f 3 {1+r}<p\le \infty,
r\in (0,1]$, and $(p,r)\neq (\infty,1)$ Then $u$ can be decomposed
as
$$
u=u^l+u^h\quad \textrm{with }u^l\in L^1(0,T; \textbf{Lip}) \textrm{
and } u^h\in L^{\widetilde{q}}(0,T; L^{\widetilde{p}})
$$
for some $\widetilde{}p,\widetilde{q}$ satisfing $\f 2
{\widetilde{q}}+\f3{\widetilde{p}}=1,\widetilde{p}>3$.

\ethm

\no{\bf Proof.}\, Fix $N\in\N$ to be determined later on. We set
\beno u^l=S_Nu,\quad u^h=u-u^l. \eeno By the definition of $S_N$ and
Lemma \ref{Berstein}, we have \ben\label{B1} \|\na u^l\|_{\infty}\le
C\sum_{j\le N-1}2^{j(1+\f3p)}\|\Delta_j u\|_p \le
C2^{2(1-\f1q)N}\|u\|_{B^{r}_{p,\infty}}. \een Due to the conditions
on $(p,q,r)$, we can choose $\widetilde{p}$ such that
$$
\widetilde{p}>\max(3,p)\quad \textrm{and} \quad \f 3 p-\f 3
{\widetilde{p}}-r<0.
$$
Thus, by Lemma \ref{Berstein} \ben\label{B-1}
\|u^h\|_{\widetilde{p}}\le \sum_{j\ge N}2^{(\f3p-\f 3
{\widetilde{p}})j}\|\Delta_ju\|_p\le C2^{(\f 3 p-\f 3
{\widetilde{p}}-r)N}\|u\|_{B^{r}_{p,\infty}}. \een Now we choose
$$
N=\Bigl[\f q{2}\log_2{(e+\|u\|_{B^{r}_{p,\infty}})}\Bigr]+1.
$$
Then by (\ref{B1}), we have \ben\label{add1} \int_0^T\|\na
u^l(t)\|_{\infty}dt\le
C\int_0^T(e+\|u(t)\|_{B^{r}_{p,\infty}})^{q}dt<+\infty. \een
On the other hand, from \eqref{B-1} we get that
 \ben\label{add2}
\int_0^T\|u^h(t)\|_{\widetilde{p}}^{\widetilde{q}}dt\le
C\int_0^T(e+\|u(t)\|_{B^{r}_{p,\infty}})^{q}dt<+\infty. \een

Hence, we complete  the proof of Lemma \ref{decom} by  \eqref{add1} and \eqref{add2}.\ef

\bthm{Lemma}\label{Integral-equ} Let $u,v$ be as in Theorem
\ref{thm1}. Set $w=u-v$. Then for any $t\in [0,T]$, there holds
\beno \langle u(t),v(t)\rangle+2\int_0^t\langle\na u,\na v\rangle
dt'=\langle u_0,v_0\rangle+\int_0^t\langle w\cdot\na u,w\rangle dt'.
\eeno \ethm

\no{\bf Proof.} Lemma \ref{decom} ensures that the trilinear form
\beno F(u,v,w):=\int_0^T\int_{\R^3}u\cdot \na
w\cdot v dxdt \eeno is continuous from $(\cL_T)^2\times L^{q}(0,T;
B^{r}_{p,\infty})$ to $\R$.
Then the lemma can be proved by following the argument of Lemma 4.4 in \cite{Ger}. Here we omit the
details.\ef

Now we are in position to prove Theorem \ref{thm1}.  Since $u$ and
$v$ are the Leray-Hopf weak solutions, there hold
 \beno
&&\|u(t)\|_2^2+2\int_0^t\|\na u(t')\|_2^2dt'\le \|u_0\|_2^2,\\
&&\|v(t)\|_2^2+2\int_0^t\|\na v(t')\|_2^2dt'\le \|v_0\|_2^2. \eeno
On the other hand, Lemma \ref{Integral-equ} yields that
 \beno \langle
u(t),v(t)\rangle+2\int_0^t\langle \na u,\na v\rangle dt'=\langle
u_0,v_0\rangle+\int_0^t\langle w\cdot\na u,w\rangle dt'. \eeno
Combining the above inequalities, we obtain \ben\label{diff-energy}
&&\|w(t)\|_2^2+2\int_0^t\|\na w(t')\|_2^2dt'\nonumber\\
&&=\|u(t)\|_2^2+\|v(t)\|_2^2-2\langle u,v\rangle (t)+2\int_0^t\|\na u(t')\|_2^2dt'+2\int_0^t\|\na v(t')\|_2^2dt'\nonumber\\
&&\quad-4\int_0^t\langle \na u,\na v\rangle (t')dt'\nonumber\\
&&\le\|u_0-v_0\|^2_2-2\int_0^t\langle w\cdot\na u,w\rangle dt'. \een
We decompose $u=u^l+u^h$ as in Lemma \ref{decom} and rewrite \beno
\int_0^t\langle w\cdot\na u,w\rangle dt'=\int_0^t\langle w\cdot\na
u^l,w\rangle dt'+\int_0^t\langle w\cdot\na u^h,w\rangle dt'. \eeno
We get by H\"{o}lder inequality that \ben\label{lowestimate}
\Bigl|\int_0^t\langle w\cdot\na u^l,w\rangle
dt'\Bigr|\le\int_0^t\|w(t')\|_2^2\|\na u^l(t')\|_\infty dt'. \een
Integration by parts, we get \beno \int_0^t\langle w\cdot\na
u^h,w\rangle dt'=-\int_0^t\langle w\cdot\na w,u^h\rangle dt', \eeno
from which and the Gagliardo-Nirenberg inequality, it follows that
\beno
\Bigl|\int_0^t\langle w\cdot\na u^h,w\rangle dt'\Bigr|&\le& \int_0^t\|\na w\|_2\|w\|_{\f{2\widetilde{p}} {\widetilde{p}-2}}\|u^h\|_{\widetilde{p}}dt'\\
&\le& C\int_0^t\|\na w\|_2\|w\|_{2}^{1-\f 3 {\widetilde{p}}}\|\na w\|_2^{\f 3 {\widetilde{p}}}\|u^h\|_{\widetilde{p}}dt'\\
&\le&
C\Bigl(\int_0^t\|w(t')\|_{2}^{2}\|u^h(t')\|_{\widetilde{p}}^{\widetilde{q}}dt'\Bigr)^{\f1{\widetilde{q}}}
\Bigl(\int_0^t\|\na w(t')\|_{2}^{2}dt'\Bigr)^{1-\f1{\widetilde{q}}}\\
&\le&
C\int_0^t\|w(t')\|_{2}^{2}\|u^h(t')\|_{\widetilde{p}}^{\widetilde{q}}dt'+\int_0^t\|\na
w(t')\|_{2}^{2}dt'. \eeno
This  together with (\ref{diff-energy})
and (\ref{lowestimate}) gives \beno
&&\|w(t)\|_2^2+\int_0^t\|\na w(t')\|_2^2dt'\\
&&\le \|u_0-v_0\|^2+C\int_0^t\|w(t')\|_2^2(\|\na
u^l(t')\|_\infty+\|u^h(t')\|_{\widetilde{p}}^{\widetilde{q}}) dt'.
\eeno This jointed with  the Gronwall  inequality produces that \beno
&&\|w(t)\|_2^2+\int_0^t\|\na w(t')\|_2^2dt'\\
&&\le \|u_0-v_0\|^2\exp\{C\int_0^t(\|\na u^l(t')\|_\infty+\|u^h(t')\|_{\widetilde{p}}^{\widetilde{q}}) dt'\}\\
&&\le
\|u_0-v_0\|_2^2\exp\{C\int_0^t(e+\|u(t')\|_{B^{r}_{p,\infty}})^qdt'\}.
\eeno

This finishes the proof of Theorem \ref{thm1}. \ef

\subsection{The proof of Theorem \ref{thm2}}
Assume that $u$ and $v$ are two weak solutions of (\ref{NSequ}) on
$(0,T)$ with the  initial data $u_0$. Let $w=u-v$, $w$ satisfies the
equation in the sense of distribution \ben\label{Diffequ} w_t-\Delta
w+w\cdot\nabla u+v\cdot\nabla w+\nabla \widetilde{p}=0, \een for
some pressure $\widetilde{p}$. In what follows, we denote
$$u_j=\Delta_ju,\quad \widetilde{p}_j=\Delta_j\widetilde{p}.$$
We get by taking the operation $\Delta_j$ on both sides of
\eqref{Diffequ} that \ben\label{localdiffequ} \p_tw_j-\Delta
w_j+\Delta_j(w\cdot\nabla u)+\Delta_j(v\cdot\nabla w)+\nabla
\widetilde{p}_j=0. \een Multiplying (\ref{localdiffequ}) by $w_j$,
we get by Lemma 2.1 for $j\ge -1$ that \ben\label{localenergy}
&&\frac{1}{2}\frac{d}{dt}\|w_j(t)\|_2^2+ca_j2^{2j}\|
w_j(t)\|_2^2\nonumber\\&& \le -\big<\Delta_j(w\cdot\na
u),w_j\big>-\big<\Delta_j(v\cdot\na w)-v\cdot\na w_j,w_j\big>, \een
with $a_{-1}=0$ and $a_j=1$ for $j\ge 0$. Here we used the fact that
$$
\big<\Delta_j(v\cdot\na w),w_j\big>=\big<\Delta_j(v\cdot\na
w)-v\cdot\na w_j,w_j\big>.
$$

\no\textbf{Case 1}. $u$ and $v$ satisfy the assumption
(a).\vspace{0.2cm}

Due to $r_1+r_2>-1$, one of $r_1$ and $r_2$ must be bigger than $-\f12$.
Without loss of generality, we assume that $r_1>-\f12$.\vspace{0.1cm}

\textbf{Step 1.} Estimate of $\big<\Delta_j(w\cdot\na
u),w_j\big>$.\vspace{0.1cm}

Using the Bony's decomposition (\ref{Bonydecom}), we have \beno
\Delta_j(w\cdot\na
u)=\Delta_j(T_{w^i}\p_iu)+\Delta_j(T_{\p_iu}w^i)+\Delta_jR(w^i,\p_iu).
\eeno Considering the support of the Fourier transform of the term
$T_{w^i}\p_iu$, we have \ben\label{LHintera}
\Delta_j(T_{w^i}\p_iu)=\sum_{|j'-j|\le4}\Delta_j(S_{j'-1}w^i\p_i\Delta_{j'}u).
\een
This gives by Lemma \ref{Berstein} that
\ben\label{LHestimate-a1} \|\Delta_j(T_{w^i}\p_iu)\|_2&\lesssim&
\sum_{|j'-j|\le4}2^{j'}\sum_{k\le j'-2}\|\Delta_{k}w\|_{\f{2p_1}
{p_1-2}}
\|\Delta_{j'}u\|_{p_1}\nonumber\\
&\lesssim& \sum_{|j'-j|\le4}2^{j'}\sum_{k\le
j'-2}2^{k\f3{p_1}}\|\Delta_{k}w\|_{2}
\|\Delta_{j'}u\|_{p_1}\nonumber\\
&\lesssim& 2^{j(1-r_1)}\|u\|_{B^{r_1}_{p_1,\infty}}\sum_{j'\le
j+2}2^{j'\f3{p_1}}\|\Delta_{j'}w\|_2. \een Similarly, we have
\ben\label{HLintera}
\Delta_j(T_{\p_iu}w^i)=\sum_{|j'-j|\le4}\Delta_j(S_{j'-1}(\p_iu)\Delta_{j'}w^i),
\een Applying Lemma \ref{Berstein} to  \eqref{HLintera} yields that
 \ben\label{HLestimate-a1}
\|\Delta_j(T_{\p_iu}w^i)\|_2&\lesssim& \sum_{|j'-j|\le4}\sum_{k\le
j'-2}2^{k}\|\Delta_{k}u\|_\infty
\|\Delta_{j'}w\|_2\nonumber\\
&\lesssim&
2^{j(1-r_1+\f3{p_1})}\|u\|_{B^{r_1}_{p_1,\infty}}\sum_{|j'-j|\le
4}\|\Delta_{j'}w\|_2. \een Since $\textrm{div}w=0$, we have
\ben\label{HHintera} \Delta_jR(w^i,\p_iu)=\sum_{j',j''\ge
j-3;|j'-j''|\le 1}\pa_i\Delta_j(\Delta_{j'}w^i\Delta_{j''}u), \een
from which and Lemma \ref{Berstein}, it follows that
\ben\label{HHestimate-a1} \|\Delta_jR(w^i,\p_iu)\|_{\f{2p_1}
{p_1+2}}&\lesssim&
\sum_{j',j''\ge j-3;|j'-j''|\le 1}2^j\|\Delta_{j'}w\|_2\|\Delta_{j''}u\|_{p_1}\nonumber\\
&\lesssim& 2^j\|u\|_{B^{r_1}_{p_1,\infty}}\sum_{j'\ge
j-3}2^{-j'r_1}\|\Delta_{j'}w\|_2. \een

Summing up (\ref{LHestimate-a1})-(\ref{HHestimate-a1}), we obtain
\ben\label{rightestimate-a1} |\big<\Delta_j(w\cdot\na u),w_j\big>|
&\lesssim&2^{j(1-r_1)}\|u\|_{B^{r_1}_{p_1,\infty}} \sum_{j'\le
j}2^{j'\f3{p_1}}\|\Delta_{j'}w\|_2\|w_j\|_2\nonumber\\&&+
2^{j(1+\f3{p_1})}\|u\|_{B^{r_1}_{p_1,\infty}}\sum_{j'\ge
j}2^{-j'r_1}\|\Delta_{j'}w\|_2\|w_j\|_2. \een

\textbf{Step 2}. Estimate of $\big<\Delta_j(v\cdot\na w)-v\cdot\na
w_j,w_j\big>$.\vspace{0.1cm}

Using the Bony's decomposition (\ref{Bonydecom}), we write \beno
&&\Delta_j(v\cdot\na w)=\Delta_j(T_{v^i}\p_iw)+\Delta_j(T_{\p_iw}v^i)+\Delta_jR(v^i,\p_iw),\\
&&v\cdot\na w_j=T_{v^i}\p_iw_j+T_{\p_iw_j}'v^i. \eeno Then we have
\beno
&&\Delta_j(v\cdot\na w)-v\cdot\na w_j\\
&&=[\Delta_j,T_{v^i}]\p_iw+\Delta_j(T_{\p_iw}v^i)+\Delta_jR(v^i,\p_iw)-T_{\p_iw_j}'v^i.
\eeno

Similar arguments as in deriving  (\ref{LHestimate-a1}) and
(\ref{HHestimate-a1}), we have \ben\label{LHestimate-a2}
&&\|\Delta_j(T_{\p_iw}v^i)\|_2
\lesssim 2^{-jr_2}\|v\|_{B^{r_2}_{p_2,\infty}}\sum_{j'\le j+2}2^{j'(1+\f3{p_2})}\|\Delta_{j'}w\|_2,\\
&&\|\Delta_jR(v^i,\p_iw)\|_{\f{2p_2} {p_2+2}} \lesssim
2^j\|v\|_{B^{r_2}_{p_2,\infty}}\sum_{j'\ge
j-3}2^{-j'r_2}\|\Delta_{j'}w\|_2.\label{HHestimate-a2} \een In view
of  the definition of $T'_{\pa_i w_j}v^i$,
\begin{align}
T'_{\pa_i w_j}v^i=\sum_{j'\ge j-2} S_{j'+2}\Delta_j\pa_i
w\Delta_{j'}v^i,\nonumber
\end{align}
and  note that $S_{j'+2}\Delta_j w=\Delta_jw$ for $j'>j$, we get
$$
\langle T'_{\pa_i w_j}v^i,w_j\rangle=\sum_{j-2\le j'\le j}\langle
S_{j'+2}\Delta_j\pa_i w\Delta_{j'}v^i,w_j\rangle,
$$
from which and Lemma \ref{Berstein}, it follows that
\ben\label{HHestimate-a22} |\langle T'_{\pa_i
w_j}v^i,w_j\rangle|\lesssim
2^{j(1+\f3{p_2}-r_2)}\|v\|_{B^{r_2}_{p_2,\infty}}\|w_j\|_2^2. \een

Now, we turn to estimate $[T_{v^i}, \Delta_j]\pa_i w$.  In view of
the definition of $\Delta_j$, we write
\begin{align}\label{commtator-representation}
[T_{v^i}, \Delta_j]\pa_i w&=\sum_{|j'-j|\le4}[S_{j'-1}v^i, \Delta_j]\pa_i w_{j'}\nonumber\\
&=\sum_{|j'-j|\le4}2^{3j}\int_{\R^3}h(2^j(x-y))(S_{j'-1}v^i(x)-S_{j'-1}v^i(y))\pa_i w_{j'}(y)dy\nonumber\\
&=\sum_{|j'-j|\le4}2^{4j}\int_{\R^3}\int_0^1y\cdot\na
S_{j'-1}v^i(x-\tau y)d\tau\pa_i h(2^jy)w_{j'}(x-y)dy,
\end{align}
from which and the Minkowski inequality, we deduce that
\ben\label{HLestimate-a2} \|[T_{v^i}, \Delta_j]\pa_i w\|_2
&\lesssim& \sum_{|j'-j|\le4}\|\na
S_{j'-1}v\|_\infty\|w_{j'}\|_2\nonumber\\
&\lesssim& 2^{j(1+\f
3{p_2}-r_2)}\|v\|_{B^{r_2}_{p_2,\infty}}\sum_{|j'-j|\le4}\|w_{j'}\|_2.
\een

Summing up (\ref{LHestimate-a2})-(\ref{HLestimate-a2}), we obtain
\ben\label{rightestimate-a2}
&&|\big<\Delta_j(v\cdot\na w)-v\cdot\na w_j,w_j\big>|\nonumber\\
&&\qquad\lesssim 2^{-jr_2}\|v\|_{B^{r_2}_{p_2,\infty}} \sum_{j'\le
j}2^{j'(1+\f3{p_2})}\|\Delta_{j'}w\|_2\|w_j\|_2\nonumber\\&&\qquad\qquad+
2^{j(1+\f3{p_2})}\|v\|_{B^{r_2}_{p_2,\infty}}\sum_{j'\ge
j}2^{-j'r_2}\|\Delta_{j'}w\|_2\|w_j\|_2. \een

Under the assumption (a), we can choose $s$ such that
\ben\label{s-condition} -r_1<s<\min(1+r_1,1+r_2). \een From
(\ref{localenergy}), (\ref{rightestimate-a1}) and
(\ref{rightestimate-a2}), it follows that
\ben\label{energyestimate-a}
&&2^{-2js}\|w_j(t)\|_2^2+a_j2^{2j(1-s)}\int_0^t\|w_j(t')\|_2^2dt'\nonumber\\
&&\lesssim\int_0^t\|u\|_{B^{r_1}_{p_1,\infty}}2^{j(1-r_1-2s)}
\sum_{j'\le j}2^{j'\f3{p_1}}\|\Delta_{j'}w\|_2\|w_j\|_2dt'\nonumber\\
&&\quad+\int_0^t\|u\|_{B^{r_1}_{p_1,\infty}}
2^{j(1+\f3{p_1}-2s)}\sum_{j'\ge j}2^{-j'r_1}\|\Delta_{j'}w\|_2\|w_j\|_2dt'\nonumber\\
&&\quad+\int_0^t\|v\|_{B^{r_2}_{p_2,\infty}}2^{-j(r_2+2s)}
\sum_{j'\le j}2^{j'(1+\f3{p_2})}\|\Delta_{j'}w\|_2\|w_j\|_2dt'\nonumber\\
&&\quad+\int_0^t\|v\|_{B^{r_2}_{p_2,\infty}}
2^{j(1+\f3{p_2}-2s)}\sum_{j'\ge j}2^{-j'r_2}\|\Delta_{j'}w\|_2\|w_j\|_2dt'\nonumber\\
&&:=I+II+III+IV. \een We set
$$W(t)=\sup_{j\ge -1}2^{-js}\|w_j(t)\|_2.$$
Using (\ref{s-condition}) and the Young's inequality, we obtain
\beno I&\le& \sum_{j'\le
j}2^{(j'-j)(r_1+s)}\int_0^t\|u\|_{B^{r_1}_{p_1,\infty}}W(t')^{\f
{2}{q_1}}
\bigl(2^{j'(1-s)}\|\Delta_{j'}w\|_2\bigr)^{1-\f 2{q_1}}2^{j(1-s)}\|w_j\|_2dt'\\
&\le&
C\Bigl(\int_0^t\|u\|_{B^{r_1}_{p_1,\infty}}^{q_1}W(t')^2dt'\Bigr)^{\f
{1}{q_1}}
\Bigl(\sup_{j\ge -1}2^{2j(1-s)}\int_0^t\|w_j(t')\|^2_2dt'\Bigr)^{\f 1 {q_1'}}\\
&\le&
C\int_0^t\|u\|_{B^{r_1}_{p_1,\infty}}^{q_1}W(t')^2dt'+\delta\sup_{j\ge
-1}2^{2j(1-s)}\int_0^t\|w_j(t')\|^2_2dt', \eeno and for $II$, we have
\beno II&\le& \sum_{j'\ge
j}2^{(j'-j)(s-1-r_1)}\int_0^t\|u\|_{B^{r_1}_{p_1,\infty}}W(t')^{\f
{2}{q_1}}
2^{j'(1-s)}\|\Delta_{j'}w\|_2\bigl(2^{j(1-s)}\|w_j\|_2\bigr)^{1-\f 2{q_1}}dt'\\
&\le&
C\int_0^t\|u\|_{B^{r_1}_{p_1,\infty}}^{q_1}W(t')^2dt'+\delta\sup_{j\ge
-1}2^{2j(1-s)}\int_0^t\|w_j(t')\|^2_2dt', \eeno and similarly for
$IV$, \beno IV\le
C\int_0^t\|v\|_{B^{r_2}_{p_2,\infty}}^{q_2}W(t')^2dt'+\delta\sup_{j\ge
-1}2^{2j(1-s)}\int_0^t\|w_j(t')\|^2_2dt', \eeno and for $III$, \beno
III&\le& \sum_{j'\le
j}2^{(j'-j)(1+r_2+s)}\int_0^t\|v\|_{B^{r_2}_{p_2,\infty}}W(t')^{\f
{2}{q_2}}
\bigl(2^{j'(1-s)}\|\Delta_{j'}w\|_2\bigr)^{1-\f 2{q_2}}2^{j(1-s)}\|w_j\|_2dt'\\
&\le&
C\int_0^t\|v\|_{B^{r_2}_{p_2,\infty}}^{q_2}W(t')^2dt'+\delta\sup_{j\ge
-1}2^{2j(1-s)}\int_0^t\|w_j(t')\|^2_2dt'. \eeno
Collecting these estimates with
(\ref{energyestimate-a}) implies that \beno W(t)^2\le
C\int_0^t(\|u(t')\|_{B^{r_1}_{p_1,\infty}}^{q_1}+\|v(t')\|_{B^{r_2}_{p_2,\infty}}^{q_1})W(t')^2dt'.
\eeno
This together with the Gronwall inequality shows that
$$
W(t)=0, \quad i.e.\quad u=v=0.
$$

This completes the proof of case (a).\vspace{0.2cm}

\no\textbf{Case 2}. $u$ and $v$ satisfy the assumption (b).
\vspace{0.1cm}

Since $u$ and $v$ are non-Lipshcitz vectors, we will use the idea of
the losing derivative estimate which was firstly introduced by
Chemin and Lerner\cite{CL}. We can refer to \cite{Dan} for a
systematic study. Recently, Danchin and Paicu\cite{DP} applied this
idea to prove the uniqueness of weak solution for the 2-D Boussinesq
equations with partial viscosity. The present proof is motivated by
\cite{DP}. We also refer to \cite{Can,Lin, Mas} for the other applications about the losing derivative estimate.

Let $s\in (0,1)$. For $\lambda>0$, we set \beno
W_j^\lambda(t)=2^{-js}e^{-\lambda\varepsilon_j(t)}\|w_j(t)\|_2,
\eeno where $\varepsilon_j(t)$ is defined by \beno
\varepsilon_j(t)=\int_0^t2^{j'}\sum_{j'\le
j+4}\bigl(\|\Delta_{j'}u(t')\|_\infty+\|\Delta_{j'}v(t')\|_\infty\bigr)dt'.
\eeno We get by (\ref{localenergy}) that
\ben\label{derilosingenergy} &&\f d
{dt}W_j^\lambda(t)+\lambda\varepsilon_j'(t)W_j^\lambda(t)+a_j2^{2j}W_j^\lambda(t)
\lesssim
2^{-js}e^{-\lambda\varepsilon_j(t)}\bigl(\|\Delta_j(w\cdot\na
u)\|_2\nonumber\\&& \qquad+\|\Delta_j(v\cdot\na w)-v\cdot\na
w_j+\sum_{j'>j}\pa_i w_j\Delta_{j'}v^i\|_2\bigr). \een Here we used
the fact that \beno
\langle\pa_iw_j\Delta_{j'}v^i,w_j\rangle=-\langle\Delta_{j'}\pa_iv^iw_j,w_j\rangle=0.
\eeno Since $W_j^\lambda(0)=0$, we get by integrating
(\ref{derilosingenergy}) on $[0,t]$ that
\ben\label{integrallosingenergy}
&&W_j^\lambda(t)+\lambda\int_0^t\varepsilon_j'(t')W_j^\lambda(t')dt'+a_j2^{2j}\int_0^tW_j^\lambda(t')dt'\nonumber\\
&&\lesssim 2^{-js}\int_0^te^{-\lambda\varepsilon_j(t')}
\|\Delta_j(w\cdot\na u)(t')\|_2dt'\nonumber\\&&
\quad+2^{-js}\int_0^te^{-\lambda\varepsilon_j(t')}\|\Delta_j(v\cdot\na
w)-v\cdot\na w_j+\sum_{j'>j}\pa_i w_j\Delta_{j'}v^i\|_2(t')dt'. \een

\no{\bf Step 1.} Estimate of $\|\Delta_j(w\cdot\na
u)\|_2$.\vspace{0.2cm}

Using the Bony's decomposition (\ref{Bonydecom}), we write \beno
\Delta_j(w\cdot\na
u)=\Delta_j(T_{w^i}\p_iu)+\Delta_j(T_{\p_iu}w^i)+\Delta_jR(w^i,\p_iu).
\eeno By (\ref{LHintera}) and Lemma \ref{Berstein}, we get
\ben\label{LHestimate-b1} \|\Delta_j(T_{w^i}\p_iu)\|_2&\lesssim&
\sum_{|j'-j|\le4}2^{j'}\sum_{k\le j'-2}\|\Delta_{k}w\|_2
\|\Delta_{j'}u\|_\infty\nonumber\\
&\lesssim&\sum_{|j'-j|\le4}2^{j'}\sum_{k\le
j'-2}2^{ks}e^{\lambda\varepsilon_k(t)}
W_k^\lambda(t)\|\Delta_{j'}u\|_\infty\nonumber\\
&\lesssim& \sum_{j'\le
j+2}2^{j's}e^{\lambda\varepsilon_{j'}(t)}W_{j'}^\lambda(t)\varepsilon_j'(t).
\een By (\ref{HLintera}), (\ref{HHintera}) and Lemma \ref{Berstein},
we have \ben\label{HLestimate-b1}
\|\Delta_j(T_{\p_iu}w^i)\|_2&\lesssim& \sum_{|j'-j|\le4}\sum_{k\le
j'-2}2^{k}\|\Delta_{k}u\|_\infty
\|\Delta_{j'}w\|_2\nonumber\\
&\lesssim&\sum_{|j'-j|\le4}2^{j's}e^{\lambda\varepsilon_{j'}(t)}
W_{j'}^\lambda(t)\sum_{k\le j'-2}2^k\|\Delta_{k}u\|_\infty\nonumber\\
&\lesssim&
\sum_{|j'-j|\le4}2^{j's}e^{\lambda\varepsilon_{j'}(t)}W_{j'}^\lambda(t)\varepsilon_j'(t),
\een and \ben\label{HHestimate-b1}
\|\Delta_jR(w^i,\p_iu)\|_2&\lesssim&
\sum_{j',j''\ge j-3;|j'-j''|\le 1}2^j\|\Delta_{j'}w\|_2\|\Delta_{j''}u\|_\infty\nonumber\\
&\lesssim& \sum_{j',j''\ge j-3;|j'-j''|\le
1}2^{j's+j}e^{\lambda\varepsilon_{j'}(t)}W_{j'}^\lambda(t)
\|\Delta_{j''}u\|_\infty\nonumber\\
&\lesssim& \sum_{j'\ge
j-3}2^{j'(s-1)+j}e^{\lambda\varepsilon_{j'}(t)}W_{j'}^\lambda(t)\varepsilon_{j'}'(t).
\een

Summing up (\ref{LHestimate-b1})-(\ref{HHestimate-b1}), we obtain
\ben\label{losing-b1} &&2^{-js}\int_0^te^{-\lambda\varepsilon_j(t')}
\|\Delta_j(w\cdot\na u)(t')\|_2dt'\nonumber\\
&&\quad\lesssim \sum_{j'\le j}2^{(j'-j)s}
\int_0^te^{\lambda(\varepsilon_{j'}(t')-\varepsilon_{j}(t'))}W_{j'}^\lambda(t')\varepsilon_j'(t')dt'\nonumber\\
&&\qquad+\sum_{j'\ge j}2^{-(j'-j)(1-s)}
\int_0^te^{\lambda(\varepsilon_{j'}(t')-\varepsilon_{j}(t'))}W_{j'}^\lambda(t')\varepsilon_j'(t')dt'.
\een

{\bf Step 2.} Estimate of $\|\Delta_j(v\cdot\na w)-v\cdot\na
w_j+\disp\sum_{j'>j}\pa_i w_j\Delta_{j'}v^i\|_2$.\vspace{0.1cm}

Using the Bony's decomposition (\ref{Bonydecom}), we write \beno
&&\Delta_j(v\cdot\na w)-v\cdot\na w_j\\
&&=[\Delta_j,T_{v^i}]\p_iw+\Delta_j(T_{\p_iw}v^i)+\Delta_jR(v^i,\p_iw)-T_{\p_iw_j}'v^i.
\eeno Similar to the proof of (\ref{LHestimate-b1}) and
(\ref{HHestimate-b1}), we get \ben\label{LHestimate-b2}
&&\|\Delta_j(T_{\p_iw}v^i)\|_2
\lesssim \sum_{j'\le j+2}2^{j's}e^{\lambda\varepsilon_{j'}(t)}W_{j'}^\lambda(t)\varepsilon_j'(t),\\
&&\|\Delta_jR(v^i,\p_iw)\|_2 \lesssim \sum_{j'\ge
j-3}2^{j'(s-1)+j}e^{\lambda\varepsilon_{j'}(t)}W_{j'}^\lambda(t)\varepsilon_{j'}'(t).\label{HHestimate-b2}
\een
Using the formula (\ref{commtator-representation}) again, we have
\ben\label{HLestimate-b2} \|[\Delta_j,T_{v^i}]\p_iw\|_2 \lesssim
\sum_{|j'-j|\le4}2^{j's}e^{\lambda\varepsilon_{j'}(t)}W_{j'}^\lambda(t)\varepsilon_j'(t).
\een Note that
$$
T_{\p_iw_j}'v^i-\sum_{j'>j}\pa_i w_j\Delta_{j'}v^i=\sum_{j-2\le
j'\le j}S_{j'+2}\Delta_j\pa_i w\Delta_{j'}v^i,
$$
it gives by Lemma \ref{Berstein} that
\ben\label{HHestimate-b22} \|T_{\p_iw_j}'v^i-\sum_{j'>j}\pa_i
w_j\Delta_{j'}v^i\|_2\lesssim
2^{js}e^{\lambda\varepsilon_{j}(t)}W_{j}^\lambda(t)\varepsilon_{j}'(t).
\een

Summing up (\ref{LHestimate-b2})-(\ref{HHestimate-b22}), we obtain
\ben\label{losing-b2} &&2^{-js}\int_0^te^{-\lambda\varepsilon_j(t')}
\|\Delta_j(v\cdot\na w)-v\cdot\na w_j)(t')\|_2dt'\nonumber\\
&&\quad\lesssim \sum_{j'\le j}2^{(j'-j)s}
\int_0^te^{\lambda(\varepsilon_{j'}(t')-\varepsilon_{j}(t'))}W_{j'}^\lambda(t')\varepsilon_j'(t')dt'\nonumber\\
&&\qquad+\sum_{j'\ge j}2^{-(j'-j)(1-s)}
\int_0^te^{\lambda(\varepsilon_{j'}(t')-\varepsilon_{j}(t'))}W_{j'}^\lambda(t')\varepsilon_{j'}'(t')dt'.
\een

From (\ref{integrallosingenergy}), (\ref{losing-b1}) and
(\ref{losing-b2}), it follows that \ben\label{losingenergy}
&&W_j^\lambda(t)+\lambda\int_0^t\varepsilon_j'(t')W_j^\lambda(t')dt'+a_j2^{2j}\int_0^tW_j^\lambda(t')dt'\nonumber\\
&&\quad\lesssim \sum_{j'\le j}2^{(j'-j)s}
\int_0^te^{\lambda(\varepsilon_{j'}(t')-\varepsilon_{j}(t'))}W_{j'}^\lambda(t')\varepsilon_j'(t')dt'\nonumber\\
&&\qquad+\sum_{j'\ge j}2^{-(j'-j)(1-s)}
\int_0^te^{\lambda(\varepsilon_{j'}(t')-\varepsilon_{j}(t'))}W_{j'}^\lambda(t')\varepsilon_{j'}'(t')dt'\nonumber\\
&&\quad :=I+II. \een
Write
$$\varepsilon_j'(t')=\varepsilon_{j'}'(t')+(\varepsilon_j'(t')-\varepsilon_{j'}'(t')),$$
and note that $\varepsilon_j'(t')-\varepsilon_{j'}'(t')\ge 0$ for
$j\ge j'$, we obtain \ben\label{I-estimate} I\lesssim \sum_{j'\le
j}2^{(j'-j)s} \int_0^tW_{j'}^\lambda(t')\varepsilon_{j'}'(t')dt' +\f
1{\lambda}\sum_{j'\le j}2^{(j'-j)s}\sup_{t'\in
[0,t]}W_{j'}^\lambda(t'), \een here we used the inequality \beno
\int_0^te^{\lambda(\varepsilon_{j'}(t')-\varepsilon_{j}(t'))}(\varepsilon_j'(t')-\varepsilon_{j'}'(t'))dt'\lesssim
\f 1\lambda, \quad \textrm{for }j'\le j. \eeno Since
$\varepsilon_{j'}(t')-\varepsilon_{j}(t')$ is an increasing function
in $t'$ for $j'\ge j$, we have \ben\label{II-estimate1} II\lesssim
\sum_{j'\ge
j}2^{-(j'-j)(1-s)}e^{\lambda(\varepsilon_{j'}(t)-\varepsilon_{j}(t))}
\int_0^tW_{j'}^\lambda(t')\varepsilon_{j'}'(t')dt'. \een

Let us for the moment assume that
\ben\label{assum}
\lambda(\|u\|_{L^1(0,t;B^1_{\infty,\infty})}+\|v\|_{L^1(0,t;B^1_{\infty,\infty})})<(1-s)\log
2.
\een
Notice that
$$
\varepsilon_{j'}(t)-\varepsilon_{j}(t)\le
(j'-j)(\|u\|_{L^1(0,t;B^1_{\infty,\infty})}+\|v\|_{L^1(0,t;B^1_{\infty,\infty})}),
$$
which together with (\ref{II-estimate1}) ensures that
\ben\label{II-estimate2} II\lesssim
\int_0^tW_{j'}^\lambda(t')\varepsilon_{j'}'(t')dt'. \een
Summing up (\ref{losingenergy}), (\ref{I-estimate}) and
(\ref{II-estimate2}), we obtain \beno
&&\sup_{j\ge -1,t'\in [0,t]}W_j^\lambda(t')+\lambda\sup_{j\ge -1}\int_0^t\varepsilon_j'(t')W_j^\lambda(t')dt'+\sup_{j\ge -1}2^{2j}\int_0^tW_j^\lambda(t')dt'\\
&&\qquad\le C\sup_{j\ge
-1}\int_0^t\varepsilon_j'(t')W_j^\lambda(t')dt'+\f
C{\lambda}\sup_{j\ge -1,t'\in [0,t]}W_j^\lambda(t'), \eeno from
which, we get by taking $\lambda$ big enough that
$$
\sup_{j\ge -1,t'\in [0,t]}W_j^\lambda(t')=0.
$$

On the other hand, the assumption (b) ensures that we can choose $t>0$ small enough such that (\ref{assum}) holds.
Thus, $u=v$ on $[0,t]$, and then we can conclude that $u=v$
on $[0,T]$ by a standard continuity argument. The proof of case (b)
is completed.\ef

\section*{Acknowledgements}
Z. Zhang would like to thank Marius Paicu  for the helpful discussions
about the losing derivative estimates. This paper was written while Z. Zhang was visiting the Mathematics
Department of Paris-Sud University as a postdoctoral fellow. He would like to thank
the hospitality and support of the Department. Q. Chen and C. Miao
were partially supported by the NSF of China under grant No.10701012, No.10725102. Z. Zhang was partially supported by the NSF
of China under grant No.10601002.

\end{document}